\documentclass[11pt,a4paper]{article}
\usepackage[T1]{fontenc}
\usepackage{lmodern}
\usepackage{HongLiu}
\usepackage{flafter}
\usetikzlibrary{arrows.meta}

\theoremstyle{plain}
\newtheorem{pathconjecture}[theorem]{Conjecture}
\crefname{pathconjecture}{conjecture}{conjectures}
\Crefname{pathconjecture}{Conjecture}{Conjectures}
\numberwithin{equation}{section}
\hypersetup{
  hidelinks,
  pdftitle={An asymptotic form of Thomasse's directed path conjecture},
  pdfsubject={Long directed paths in oriented graphs of large minimum out-degree},
  pdfkeywords={oriented graph, directed path, minimum out-degree, longest cycle}}

\title{An asymptotic form of Thomass\'e's directed path conjecture}
\author{%
  Xinyao Ren,
  Yisai Xue\thanks{Email: \texttt{xueyisai@nbu.edu.cn}.  Supported by the National Natural Science Foundation of China (No. 12501486).},
  and Ai Yang\\[4pt]
  {\small School of Mathematics and Statistics, Ningbo University, Ningbo, China}%
}
\date{}

\begin{document}
\maketitle

\begin{abstract}
Thomass\'e conjectured that every oriented graph of minimum out-degree
$d$ contains a directed path of length $2d$. We resolve this conjecture
asymptotically by proving that every such graph contains a directed
path of length at least $2d-O(d\log\log d/\log d)$.
\end{abstract}

\section{Introduction}

Every finite loopless digraph of minimum out-degree $d$ has a directed path of length at least $d$: all out-neighbours of the last vertex of a longest path lie on it. Replacing each edge of $K_{d+1}$ by two opposite arcs shows sharpness. In 2005, Thomass\'e~\cite{Thomasse} conjectured that the bound doubles for oriented graphs (orientations of finite simple graphs); see also Bondy and Murty~\cite[Appendix~A, Problem~70, p.~589]{BM}.

\begin{pathconjecture}[Thomass\'e]\label[pathconjecture]{conj:path}
Let $d$ be a positive integer. Every oriented graph of minimum out-degree at least $d$ contains a directed path of length $2d$.
\end{pathconjecture}

Regular tournaments on $2d+1$ vertices show that the conjectured bound
is best possible. Jackson~\cite{Jackson} proved \Cref{conj:path} when
the minimum in-degree is also at least $d$. Without an in-degree
assumption, Thomass\'e~\cite{Thomasse} obtained the bound $3d/2$,
for which Cheng and Keevash~\cite{CK} later published a proof.
For related results on other subgraphs forced by large minimum
out-degree, see Hons et al.~\cite{HonsEtAl}. We prove the following
asymptotic form of \Cref{conj:path}.

\begin{theorem}\label{main}
There is an absolute constant $K>0$ such that, for all sufficiently large integers $d$, every oriented graph of minimum out-degree at least $d$ contains a directed path of length at least $2d-K\,\frac{d\log\log d}{\log d}$.
\end{theorem}

\paragraph{Proof idea.}
Suppose that $D$ is strongly connected, has minimum out-degree $d$,
and contains no path of length $2d$. A longest cycle $C$ is then
nonspanning, and a longest path in $D-C$ extends through $C$ by
$c=|C|$ arcs. Our deletion lemma finds a nonempty induced subgraph
of $D-C$ with minimum out-degree at least $d-c/2-o(c)$. Comparing
this degree loss with the $c$ arcs gained gives the leading factor
$2$ in the induction on $|V(D)|$. The quantitative form of the
lemma gives the stated error term.

If the deletion bound fails, we sample blocks of consecutive vertices
of $C$ and construct an auxiliary digraph $H$ on their union $X$.
An arc $x\to y$ of $H$ represents a path from $x$ to the successor
of $y$ on $C$, with all internal vertices outside $C$. In the cyclic
order on $X$, we find disjoint intervals $I_1,\ldots,I_j$ and
$j$ vertex-disjoint $j$-cycles of $H$, each visiting these intervals
in the indicated order.
Within each interval the cycles appear in a common order, so the corresponding paths rejoin all
the segments of $C$ into a single cycle. At least one of these
paths visits a vertex outside $C$, giving a longer cycle.

\subsection{Notation}

In a digraph $D$, write $u\to v$ for an arc, $N_D^+(v)$ and $N_D^-(v)$ for the out- and in-neighbourhoods of $v$, and $\delta^+(D)$ for the minimum out-degree. Let $D[S]$ denote the subdigraph induced by $S$.

All paths, cycles and walks are directed, with length measured in arcs. Let $\ell(D)$ be the maximum path length. For a cycle $C$, write $|C|=|V(C)|$ and $D-C=D[V(D)\setminus V(C)]$. For $y\in V(C)$, let $y^+$ be its successor on $C$, and for $T\subseteq V(C)$ let $T^+=\{y^+:y\in T\}$. All logarithms are natural.

\section{From cycle deletion to long paths}\label{sec:delete}

The next two lemmas compare the degree lost on deleting a longest cycle with the length gained by extending a path through it.

\begin{lemma}[Cycle deletion]\label[lemma]{corelemma}
Let $r\ge16$ be an integer and let $D$ be an oriented graph with $\delta^+(D)\ge d$. Let $C$ be a longest directed cycle of $D$. If $V(C)\ne V(D)$ and $c=|C|\ge r^{10r}$, then $D-C$ contains a nonempty induced subgraph $J$ with
\[
  \delta^+(J)\ge d-(1/2+8/r)c.
\]
\end{lemma}

The threshold $r^{10r}$ is used in the two counting estimates in \Cref{sec:proof}.

\begin{lemma}[Path extension]\label[lemma]{lem:extension}
Let $D$ be a strongly connected oriented graph with $\delta^+(D)\ge d\ge1$
and $\ell(D)<2d$. Every longest cycle $C$ is nonspanning and satisfies
\[
  |C|\ge d+1 \qquad\text{and}\qquad |C|+\ell(D-C)\le\ell(D).
\]
\end{lemma}

\begin{proof}
The last vertex of a longest path has all its out-neighbours on the
path. Its earliest out-neighbour closes a cycle of length at least
$d+1$, so $|C|\ge d+1$. Counting arcs gives $|V(D)|\ge2d+1$, so a
spanning cycle would yield a path of length at least $2d$.

\begin{claim}[{\cite[Claim~10]{CK}}]
\label[claim]{cl:long-cycles-intersect}
Any two cycles of length at least $d+1$ share a vertex.
\end{claim}
Indeed, if $C_1$ and $C_2$ were vertex-disjoint, strong connectivity
would give a shortest path from $C_1$ to $C_2$, with all internal
vertices outside both cycles. Extending this path backwards along
$C_1$ and forwards along $C_2$ gives a path of length at least
$|C_1|+|C_2|-1\ge 2d+1$, contrary to $\ell(D)<2d$.

Let $v$ be the last vertex of a longest path $P$ in $D-C$. If $v$ has
no out-neighbour on $C$, all its out-neighbours lie on $P$. The earliest
one closes a cycle of length at least $d+1$ disjoint from $C$, contrary
to Claim~\ref{cl:long-cycles-intersect}. Thus $v$ sends an arc to $C$.
Taking this arc and then $|C|-1$ arcs of $C$ extends $P$ by $|C|$ arcs.
\end{proof}

\begin{proof}[\textbf{Proof of \Cref{main}}]
Fix an integer $r\ge16$. We prove by induction on $|V(D)|$ that every nonempty oriented graph $D$ satisfies
\begin{equation}\label{tradeoff}
  \delta^+(D)\le (1/2+8/r)\ell(D)+r^{10r}.
\end{equation}
Let $d=\delta^+(D)$. The inequality holds if $d\le r^{10r}$ or $\ell(D)\ge2d$. 
If $D$ is not strongly connected, choose a \textit{sink strongly connected component} $H$ of $D$, that is, one with no arc from $H$ to any other component. 
Then $\delta^+(H)\ge d$ and $\ell(H)\le\ell(D)$, so the induction hypothesis applied to $H$ gives the desired bound.
Hence assume that $D$ is strongly connected, $d>r^{10r}$ and $\ell(D)<2d$.

Choose a longest cycle $C$ of order $c$. By \Cref{lem:extension}, $C$ is nonspanning, $c\ge d+1$ and $\ell(D-C)+c\le\ell(D)$. Apply \Cref{corelemma} to obtain $J$. Induction inside $J$ gives
\[
  d\le\delta^+(J)+(1/2+8/r)c
   \le (1/2+8/r)\bigl(\ell(J)+c\bigr)+r^{10r}
   \le (1/2+8/r)\ell(D)+r^{10r}.
\]
This proves \eqref{tradeoff}. For sufficiently large $d\le\delta^+(D)$, let $r=\lfloor\log d/(20\log\log d)\rfloor$. Then $r\ge16$ and $10r\log r\le(\log d)/2$, so $r^{10r}\le\sqrt d$. Equation~\eqref{tradeoff} gives $\ell(D)\ge 2(d-\sqrt d)/(1+16/r)$. Since $1/r=O(\log\log d/\log d)$, the theorem follows.
\end{proof}

\section{Proof of the cycle deletion lemma (Lemma \ref{corelemma})}\label{sec:proof}

Let $L=V(D)\setminus V(C)$ and $k=d-(1/2+8/r)c$. Suppose for a contradiction that every nonempty induced subgraph of $D[L]$ has minimum out-degree less than $k$. As $L$ is nonempty, $k>0$. We shall build a cycle containing all vertices of $C$ and at least one vertex of $L$.

Colour a vertex blue if it has at least $(1/2+6/r)c$ out-neighbours on $C$, and red otherwise. 
Every red vertex has more than $d-(1/2+6/r)c=k+2c/r$ out-neighbours in $L$. We shall link red vertices on $C$ to blue vertices in $L$, which have many choices for an arc back to $C$.

\subsection{Constructing the auxiliary digraph}

We choose $X\subseteq V(C)$ and construct an auxiliary digraph $H$ on $X$.
To reconnect $C$ after deleting selected arcs $x\to x^+$, we record each
chosen path in $D$ from $x$ to $y^+$ with all internal vertices in $L$
as an arc $x\to y$ of $H$.
Let $m=r\lfloor c/r^2\rfloor$. Then $c/(2r)\le m\le c/r$, so
$m\ge r^{8r}$. We first restrict the possible blue endpoints to a set of size at most $cm$, fixed before sampling $X$.
This allows us to impose the required concentration estimates
simultaneously without any bound on $|L|$.

\begin{theorem}[Menger, set form]\label{thm:menger}
For vertex sets $A,B$ in a digraph, the maximum number of vertex-disjoint
$A$--$B$ paths equals the minimum size of a set meeting all such paths.
A vertex of $A\cap B$ counts as a path of length zero.
\end{theorem}

This follows from \cite[Theorem~5.4.1(b)]{BG} by adding a source with
arcs to $A$ and a sink with arcs from $B$.

\begin{claim}\label[claim]{cl:blue-endpoints}
There is a set $Z\subseteq L$ of at most $cm$ blue vertices such that
any set $S$ of at most $m$ red vertices of $C$ can be joined to $Z$ by
$|S|$ vertex-disjoint paths, each lying in $L$ after its first vertex.
\end{claim}

\begin{poc}
Fix a red vertex $x\in V(C)$ and let $A=N_D^+(x)\cap L$.
Given $F\subseteq L$ with $|F|<m$, let $Y$ be the set of vertices
reachable from $A\setminus F$ in $D[L\setminus F]$.
Since $|A|>k+2c/r\ge m>|F|$, we have
$\varnothing\ne A\setminus F\subseteq Y$.
By the definition of $Y$, every out-neighbour of $y\in Y$ in
$L\setminus F$ also lies in $Y$. Thus each red vertex in $Y$ has
more than $k+2c/r-|F|>k$ out-neighbours in $Y$.
Since $\delta^+(D[Y])<k$, the set $Y$ contains a blue vertex.
Thus deleting fewer than $m$ vertices from $D[L]$ always leaves a path
from $A$ to a blue vertex. By Menger Theorem (\Cref{thm:menger}), there are $m$
vertex-disjoint such paths in $D[L]$. Adding the arc from $x$ to the
first vertex of each path gives $m$ paths that meet pairwise only at $x$.
Fix one such family for each red vertex of $C$, and let $Z$ be the set
of all their blue endpoints. There are at most $c$ families, so $|Z|\le cm$.

Now let $S$ be a nonempty set of at most $m$ red vertices on $C$,
and work in $D[S\cup L]$. Given $F\subseteq S\cup L$ with $|F|<|S|$,
choose $x\in S\setminus F$. The $m$ fixed paths from $x$ to $Z$ meet
pairwise only at $x$, so $F$ meets at most $|F|<m$ of them.
Thus some path from $S$ to $Z$ avoids $F$. By Menger Theorem (\Cref{thm:menger}),
there are $|S|$ vertex-disjoint paths from $S$ to $Z$.
Their starting vertices exhaust $S$, so every subsequent vertex lies in $L$.
\end{poc}

\begin{theorem}[Hoeffding's inequality,
  {\cite[Section~6]{Hoeffding}}]\label[theorem]{lem:hoeffding}
Let $w_1,\ldots,w_N\in[0,B]$, where $B>0$, and let $1\le q\le N$
be an integer. Choose a $q$-element subset $S$ of $\{1,\ldots,N\}$
uniformly at random, and let $Y=\sum_{i\in S}w_i$. Then
$\mathbb E Y=\frac{q}{N}\sum_{i=1}^N w_i$, and for every $t>0$,
\[
  \Pr\!\left(Y\le \mathbb E Y-t\right)
  \le e^{-2t^2/(qB^2)}.
\]
\end{theorem}

Partition $C$ into intervals of $r$ consecutive vertices and a remainder
of fewer than $r$ vertices. Choose $m/r$ of the full intervals uniformly
without replacement, and let $X$ be their union. 
If a vertex of $X$ has its successor outside $X$, then it is the last vertex of a selected interval.
Hence $|X|=m$ and $|X^+\setminus X|\le m/r$.

For a blue $v\in V(C)\cup Z$, let $Y_v=|N_D^+(v)\cap X^+|$. 
Except for the successors of the remainder, each vertex of $C$ belongs to
$X^+$ with probability at least $m/c$. Hence
\[
  \mathbb E Y_v\ge \frac mc\bigl((1/2+6/r)c-r\bigr)
  \ge (1/2+6/r)m-1.
\]
Each full interval $I$ has weight $|N_D^+(v)\cap I^+|\in[0,r]$,
and $Y_v$ is the sum of the $m/r$ sampled weights. Since $m\ge 2r$,
we have $\mathbb E Y_v-(1/2+5/r)m\ge m/r-1\ge m/(2r)$.
Applying Hoeffding's inequality (Theorem~\ref{lem:hoeffding}) with $q=m/r$, $B=r$, and $t=m/(2r)$ gives
\[
  \Pr\bigl(Y_v<(1/2+5/r)m\bigr)\le e^{-m/(2r^3)}.
\]
As $m\ge c/(2r)$ and $c\ge r^{10r}$, we have $(c+cm)e^{-m/(2r^3)}<1$. 
By the union bound, fix a choice of $X$ for which $Y_v\ge(1/2+5/r)m$ holds for every blue $v\in V(C)\cup Z$.

Let $R$ be the set of red vertices in $X$. Since $|R|\le m$,
Claim~\ref{cl:blue-endpoints} gives pairwise vertex-disjoint paths $P_x$,
one from each $x\in R$ to $Z$, each lying in $L$ after its first vertex.
For $x\in X\setminus R$, let $P_x$ consist of $x$ alone. In either case,
the last vertex $v_x$ of $P_x$ is blue and belongs to $X\cup Z$.

For distinct $x,y\in X$ with $y^+\ne x$, add the arc $x\to y$ to $H$
whenever $D$ contains the arc $v_x\to y^+$. The corresponding path
in $D$ follows $P_x$ from $x$ to $v_x$ and then takes this arc to $y^+$.
All its internal vertices lie in $L$, and paths with different starting
vertices have disjoint interiors.

The $Y_{v_x}$ out-neighbours of $v_x$ in $X^+$ give distinct choices
of $y\in X$. Only $x$ and its predecessor on $C$ can be excluded, so
$d_H^+(x)\ge Y_{v_x}-2\ge(1/2+5/r)m-2$.
Since $2\le m/r$, $|X|=m$, and $|X^+\setminus X|\le m/r$, we obtain
\begin{equation}\label{aux}
  \delta^+(H)\ge(1/2+4/r)m,\qquad |X\cup X^+|\le(1+1/r)m.
\end{equation}

\subsection{Finding closed walks through red vertices}

For $U\subseteq X$, let $\Gamma(U)$ consist of the vertices with at
least $m/r^2$ in-neighbours in $U$ in the digraph $H$. 
Each vertex in $X\setminus\Gamma(U)$ receives fewer than $m/r^2$ arcs
from $U$. Since there are at most $m$ such vertices, they receive at
most $m^2/r^2$ arcs from $U$ in total.
If $|U|\ge m/r$, then some $u\in U$ has at most
$m^2/(r^2|U|)\le m/r$ out-neighbours outside $\Gamma(U)$.
Since $u$ has at least $\delta^+(H)$ out-neighbours in total,
\eqref{aux} gives
\begin{equation}\label{eq:walk-neighbourhood}
  |\Gamma(U)|\ge \delta^+(H)-m/r\ge (1/2+3/r)m.
\end{equation}

\begin{claim}\label[claim]{cl:outgoing}
If $|U|\ge m/2$ and $|U\cap R|<m/r$, then $|\Gamma(U)|\ge |U|+m/r$.
\end{claim}

\begin{poc}
Let $W=(U\setminus R)\cap\Gamma(U)^+$. Since
$|\Gamma(U)^+|=|\Gamma(U)|$ and $|U\cap R|<m/r$, inclusion--exclusion
and \eqref{aux} give
\[
  |W|\ge |U\setminus R|+|\Gamma(U)|-|X\cup X^+|
      \ge |U|+|\Gamma(U)|-(1+2/r)m.
\]
As $|U|\ge m/2$, \eqref{eq:walk-neighbourhood} now gives $|W|\ge m/r$.

Since $W\subseteq U$, at most $m^2/r^2$ arcs of $H$ from $W$ end
outside $\Gamma(U)$. Consider the arcs from $W$ to $\Gamma(U)$.
Every $w\in W$ is blue, so each such arc $w\to y$ corresponds to
a distinct arc $w\to y^+$ in $D$. At most $\binom{|W|}{2}$ of these
arcs end in $W$, since $D[W]$ is oriented. The remaining heads lie in
$\Gamma(U)^+\setminus W$, giving at most $|\Gamma(U)|-|W|$ further
choices for each $w\in W$. Averaging the out-degrees in $H$ over $W$
therefore gives
\[
  \delta^+(H)
  \le \frac{m^2}{r^2|W|}+\frac{|W|-1}{2}+|\Gamma(U)|-|W|
  \le \frac{m}{r}+|\Gamma(U)|-\frac{|W|}{2}.
\]
Substituting $|W|\ge |U|+|\Gamma(U)|-(1+2/r)m$ and using \eqref{aux},
we obtain $|\Gamma(U)|-|U|\ge 2\delta^+(H)-(1+4/r)m\ge m/r$, as claimed.
\end{poc}

An \textit{ordered closed walk} of length $j$ in $H$ is a tuple $(v_1,\ldots,v_j)$ such that $v_i\to v_{i+1}$ in $H$ for
$1\le i\le j$, where $v_{j+1}=v_1$. Vertices may repeat, and distinct tuples are counted separately.

\begin{claim}\label[claim]{cl:closed-walks}
For some $3\le j\le r$, at least $(m/r^2)^j$ ordered closed walks of
length $j$ in $H$ meet $R$.
\end{claim}

\begin{poc}
Let $U_0=\{v\in X:|N_H^-(v)|\ge m/2\}$ and $U_{i+1}=\Gamma(U_i)$.
We have $|U_0|\ge m/r$, since otherwise the average in-degree would
be at most $m/2+m/r$, contrary to \eqref{aux}. By
\eqref{eq:walk-neighbourhood}, $|U_i|>m/2$ for every $i\ge1$.
Some $U_h$ with $1\le h\le r-2$ contains at least $m/r$ red vertices.
Otherwise Claim \ref{cl:outgoing} would give
$|U_{r-1}|\ge |U_1|+(r-2)m/r>m$.

Choose $y\in U_h\cap R$ and trace a walk backwards from $y$, choosing
one predecessor in each of $U_{h-1},\ldots,U_0$. If its start is
$v\in U_0$, then
\[
  |N_H^+(y)\cap N_H^-(v)|\ge \delta^+(H)+m/2-m\ge 4m/r.
\]
Any vertex $z$ in this intersection closes the walk by $y\to z\to v$.
For $j=h+2$, start the ordered walk at the chosen red vertex $y$.
The next vertex is $z$, and the remaining $h$ vertices record the
predecessor choices in reverse order. Each of the $j$ choices has
at least $m/r^2$ options, giving at least $(m/r^2)^j$ ordered walks.
\end{poc}

\subsection{Reassembling the cycle}

Fix $j$ as in Claim~\ref{cl:closed-walks}.
Label the vertices of $X$ by $\mathbb Z/m\mathbb Z$ in the cyclic order inherited from $C$. 
An \textit{interval} of $X$ is a set of consecutive labels modulo $m$. For $x\in X$, the vertex $x+1$ is the next vertex of $X$ in this order, and $x^+$ is the successor of $x$ on $C$. 
For labels $a,b\in\{0,\ldots,m-1\}$, their \textit{cyclic distance} in $X$ is $\min\{|a-b|,\,m-|a-b|\}$.

\begin{claim}\label[claim]{cl:ordered-cuts}
There are pairwise disjoint intervals $I_1,\ldots,I_j$ of $X$
and vertices $x_{i,t}\in I_i$ for $1\le i,t\le j$ with the
following properties. For each $i$, the vertices
$x_{i,1},\ldots,x_{i,j}$ are distinct and occur in this order
along $I_i$. For each $t$, the column
$x_{1,t}\to\cdots\to x_{j,t}\to x_{1,t}$ is a directed cycle
in $H$ meeting $R$.
\end{claim}

\begin{poc}
Let $s=\lfloor\sqrt m\rfloor$. 
By Claim~\ref{cl:closed-walks}, at least $m^j r^{-2j}$ tuples $(v_1,\ldots,v_j)\in X^j$ give a
closed walk $v_1\to\cdots\to v_j\to v_1$ in $H$ meeting $R$.
Call a tuple $(v_1,\ldots,v_j)$ \emph{separated} if the cyclic distance in $X$ between $v_p$ and $v_q$ is at least $s$ for every $1\le p<q\le j$.
For each pair of positions, once one entry
is fixed, at most $2s-1$ choices for the other lie at distance
less than $s$. Thus at most $\binom{j}{2}(2s-1)m^{j-1}
\le j^2s m^{j-1}$ tuples fail to be separated.

Choose $\mathbf a=(a_1,\ldots,a_j)$ uniformly from $X^j$. For each
$0\le h<s$, the shift $\mathbf a+h=(a_1+h,\ldots,a_j+h)$ is also
uniform. Hence the expected number of shifts that are separated
and give a closed walk meeting $R$ is at least $s r^{-2j}-j^2s^2/m\ge r^{2r}-r^2>j$,
since $j\le r$ and $r^{4r}\le s\le\sqrt m$.

Choose $\mathbf a$ and offsets $0\le h_1<\cdots<h_j<s$ giving
$j$ such shifts. Translation preserves cyclic distance, so $\mathbf a$
is separated and the intervals $I_i=\{a_i,\ldots,a_i+s-1\}$ are
pairwise disjoint. Set $x_{i,t}=a_i+h_t$. The increasing offsets
give the required order in each row. Each column gives a closed walk
meeting $R$, and its entries are distinct by separation, so it is
a directed cycle.
\end{poc}

Fix an array as in Claim~\ref{cl:ordered-cuts} and take row indices
modulo $j$. For each selected vertex $x_{i,t}$, replace the arc
$x_{i,t}\to x_{i,t}^+$ of $C$ by the path from $x_{i,t}$ to
$x_{i+1,t}^+$ represented by $x_{i,t}\to x_{i+1,t}$ in $H$.
Each selected vertex starts one replacement path, and each
$x_{i,t}^+$ receives the path from $x_{i-1,t}$. Since these paths
have pairwise disjoint interiors in $L$, every vertex of the
resulting subgraph has one incoming and one outgoing arc. We thus
obtain vertex-disjoint directed cycles covering $V(C)$. Each cycle
contains a selected vertex, since the retained arcs of $C$ form
disjoint paths.

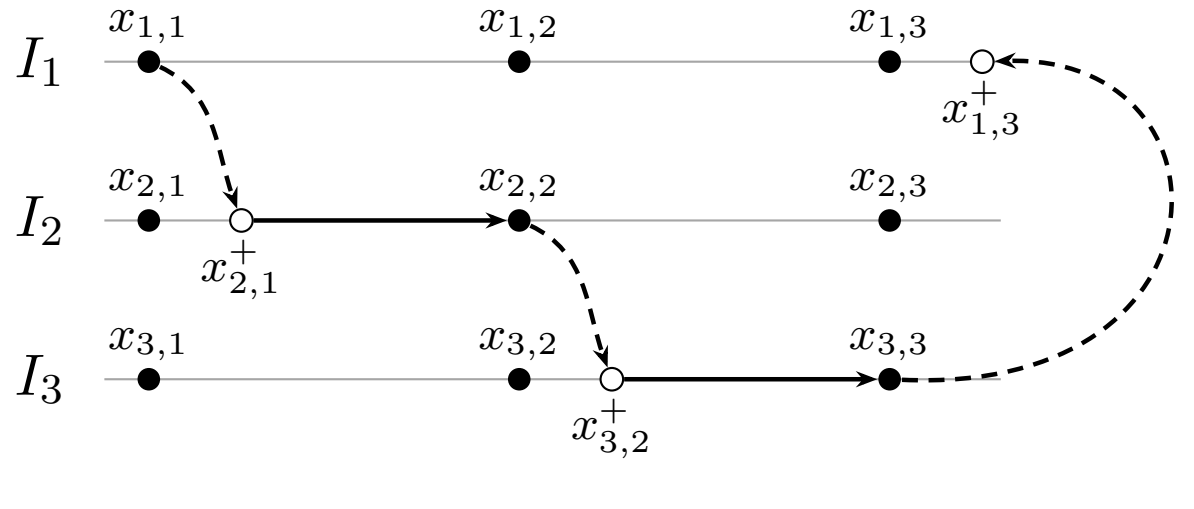
\begin{figure}[!htb]
\centering
\begin{tikzpicture}[
  x=2.45cm,y=1.05cm,
  cut/.style={circle,fill=black,inner sep=1.5pt},
  head/.style={circle,draw,fill=white,inner sep=1.5pt},
  guide/.style={draw=black!35,line width=.4pt},
  kept/.style={-{Stealth[length=4pt]},line width=.85pt},
  jump/.style={kept,densely dashed},
  every label/.style={font=\footnotesize,inner sep=1pt}]
  \foreach \row/\yy in {1/2,2/1,3/0}{
    \node[anchor=east] at (-.17,\yy) {$I_{\row}$};
    \draw[guide] (-.12,\yy)--(2.30,\yy);
    \foreach \col/\xx in {1/0,2/1,3/2}{
      \node[cut,label=above:{$x_{\row,\col}$}]
        (c\row\col) at (\xx,\yy) {};
    }
  }
  \node[head,label=below:{$x_{2,1}^+$}] (s21) at (.25,1) {};
  \node[head,label=below:{$x_{3,2}^+$}] (s32) at (1.25,0) {};
  \node[head,label=below:{$x_{1,3}^+$}] (s13) at (2.25,2) {};
  \draw[jump] (c11) to[out=-25,in=110] (s21);
  \draw[kept] (s21) -- (c22);
  \draw[jump] (c22) to[out=-25,in=110] (s32);
  \draw[kept] (s32) -- (c33);
  \draw[jump] (c33) .. controls (2.95,-.10) and (2.95,2.05) .. (s13);
\end{tikzpicture}
\caption{The passage from $x_{1,1}$ to $x_{1,3}^+$ for $j=3$.
Each row runs from left to right along $C$. Dashed arrows are replacement
paths and solid arrows are retained segments of $C$.}
\label{fig:switch}
\end{figure}

Let $F(x)$ be the next selected vertex after $x$ on its cycle in the cover.
Write $I_{\pi(i)}$ for the interval following $I_i$ along $C$.
The replacement path from $x_{i,t}$ ends at $x_{i+1,t}^+$. The retained
segment of $C$, possibly of length zero, then leads to $x_{i+1,t+1}$
if $t<j$ and to $x_{\pi(i+1),1}$ if $t=j$. Hence
\[
  F^{t-1}(x_{i,1})=x_{i+t-1,t}\quad(1\le t\le j),
  \qquad F^j(x_{i,1})=x_{\pi(i),1}.
\]
In the last step the replacement path from $x_{i+j-1,j}$ ends at
$x_{i,j}^+$, as shown in \Cref{fig:switch} for $i=1$ and $j=3$.

The first-return map on the first column acts as $\pi$, a $j$-cycle.
Thus the first column lies on one cycle $Q$. For each $t$, the vertices
$x_{i+t-1,t}$ run through column $t$ as $i$ varies, so $Q$ contains
every selected vertex. Since every cycle of the cover contains a selected
vertex, $Q$ is the only cycle and contains $V(C)$. By Claim~\ref{cl:ordered-cuts}, the first column contains a red vertex $x$. Its replacement path passes
through $v_x\in L$, so $|Q|>|C|$, a contradiction.\qed

\section*{AI declaration}

The authors used ChatGPT to refine proof strategies, check intermediate
calculations, and improve exposition. The authors independently verified all
mathematical content, reviewed and edited all AI-assisted material,
and take full responsibility for the manuscript.

\begingroup
\small
\bibliographystyle{abbrv}
\bibliography{refs}
\endgroup
\end{document}